\documentclass[reqno,12pt]{amsart}
\usepackage[margin=0.95in]{geometry}
\usepackage{amsmath,amssymb,amsthm,mathtools}
\usepackage[hidelinks]{hyperref}
\hypersetup{pdftitle={Lifts of Brauer characters and vertex pairs for odd primes},pdfauthor={Yong Yang}}

\newcommand{\Irr}{\operatorname{Irr}}
\newcommand{\IBr}{\operatorname{IBr}}
\newcommand{\Lin}{\operatorname{Lin}}
\newcommand{\Ind}{\operatorname{Ind}}
\newcommand{\Res}{\operatorname{Res}}
\newcommand{\N}{\operatorname{N}}
\newcommand{\Stab}{\operatorname{Stab}}
\newcommand{\Hom}{\operatorname{Hom}}

\newtheorem{theorem}{Theorem}[section]
\newtheorem{lemma}[theorem]{Lemma}

\newtheorem{corollary}[theorem]{Corollary}
\theoremstyle{remark}

\numberwithin{equation}{section}

\title[Lifts of Brauer characters]{Lifts of Brauer characters and vertex pairs for odd primes}
\author[Y. Yang]{Yong Yang}
\address{Department of Mathematics, Texas State University, 601 University Drive, San Marcos, TX 78666, USA}
\email{yang@txstate.edu}
\date{}

\begin{document}

\begin{abstract}
Let $p$ be an odd prime, let $G$ be a finite $p$-solvable group, and let $\varphi\in\IBr_p(G)$ have vertex $Q$. We prove that $\varphi$ has at most $|Q:Q'|$ lifts to ordinary irreducible characters, establishing Cossey's conjecture for odd primes.
\end{abstract}

\maketitle

\section{Introduction}
Let $p$ be a prime and let $G$ be a finite $p$-solvable group. If $\varphi\in\IBr_p(G)$, a \emph{lift} of $\varphi$ is a character $\chi\in\Irr(G)$ such that $\chi^0=\varphi$, where $\chi^0$ denotes the restriction of $\chi$ to the $p$-regular elements. Throughout, a vertex of an irreducible Brauer character means the Green vertex of its corresponding simple module. We write
\[
 L_\varphi=\{\chi\in\Irr(G):\chi^0=\varphi\}.
\]
For a group $Q$, let $Q'$ denote its commutator subgroup and let $\Lin(Q)$ denote the set of its linear characters. It was conjectured in \cite{Cossey2007} that if $Q$ is a vertex of $\varphi$, then
\begin{equation}\label{eq:cossey}
 |L_\varphi|\le |Q:Q'|.
\end{equation}
We prove \eqref{eq:cossey} for every odd prime. The conjecture is still open in characteristic two; one difficulty is that vertex pairs of a lift need not be linear or conjugate in that case \cite{JinWang}. For odd $p$, the vertex pairs of a lift form one conjugacy class, and their vertex characters are linear \cite[Theorem~1.2 and Lemma~4.3]{CLvertex}. Thus, if one can bound the lifts having one fixed vertex pair $(Q,\delta)$ by the size of the $\N_G(Q)$-orbit of $\delta$, then \eqref{eq:cossey} follows by summing over the orbits on $\Lin(Q)$.

The remaining step is closely related to a question of Lewis. Let $Q\le H\le G$, and let
\[
 I_\varphi(H\mid Q)=\{\eta\in\IBr_p(H):\eta^G=\varphi,\ Q\text{ is a vertex of }\eta\}.
\]
For solvable groups, the bound
\[
 |I_\varphi(H\mid Q)|\le |\N_G(Q):\N_H(Q)|
\]
was proved when either $|G|$ is odd or $p=2$; the even-order case in odd characteristic was left open \cite{LewisInducing}.

We prove a stronger module-theoretic statement that does not require solvability. Let $k$ be an algebraically closed field of characteristic $p$, let $M$ be a simple $kG$-module with Green vertex $Q$, and let $Q\le H\le G$. Then the number of isomorphism classes of simple $kH$-modules $L$ with vertex $Q$ and $\Ind_H^G L\cong M$ is at most
\[
 |N_G(Q):N_H(Q)|.
\]
For a fixed source, the Green--Puig parametrization assigns distinct indecomposable projective multiplicity modules to distinct inducing modules. We prove that these multiplicity modules induce to a common module. The required bound then follows from the Krull--Schmidt theorem and a dimension count. The $p$-solvable hypothesis is used in the passage from this inducing bound to lifts.

\section{Sources and multiplicity modules}

Let $k$ be an algebraically closed field of characteristic $p$. All modules are finite-dimensional. We use left conjugation ${}^gH=gHg^{-1}$ in Mackey decompositions. For right conjugation, we write $H^g=g^{-1}Hg$ and $E^g={}^{g^{-1}}E$.

Let $Q$ be a $p$-group and let $E$ be an indecomposable $kQ$-module. For a $kQ$-module $N$, define
\[
 \mathcal M_E(N)=\Hom_{kQ}(E,N)/\operatorname{rad}_{kQ}(E,N),
\]
where the denominator is the categorical radical. Explicitly, it consists of the homomorphisms $f:E\to N$ such that $1_E-gf$ is invertible for every $g:N\to E$. We write $J(A)$ for the Jacobson radical of an algebra $A$. If
\[
 N\cong E^m\oplus N_0
\]
and $N_0$ has no direct summand isomorphic to $E$, then
\[
 \mathcal M_E(N)\cong
 \bigl(\operatorname{End}_{kQ}(E)/J(\operatorname{End}_{kQ}(E))\bigr)^m
 \cong k^m.
\]
In particular, $\mathcal M_E$ is additive, and $\mathcal M_E(N)\ne0$ precisely when $E$ is a direct summand of $N$. We also write $\mathcal M_E(N)$ for $\mathcal M_E(N\!\downarrow_Q)$ when $N$ is a module for a group containing $Q$.

Let $G$ be a finite group containing $Q$, suppose that $E$ has vertex $Q$, and set
\[
 C=\Stab_{\N_G(Q)}([E]),
\]
where $[E]$ is the isomorphism class of $E$. We describe a choice of twisted group algebra that depends only on this source and is compatible with restriction to subgroups of $C$.

Let $\rho_E$ denote the action of $Q$ on $E$. Choose a transversal $\mathcal T$ for $Q$ in $C$ containing $1$. For each $t\in\mathcal T$, choose $u_t\in\operatorname{GL}_k(E)$ such that
\[
 u_t\rho_E(q)u_t^{-1}=\rho_E(tqt^{-1})
 \qquad(q\in Q),
\]
with $u_1=1_E$. For $c=qt$, where $q\in Q$ and $t\in\mathcal T$, define $u_c=\rho_E(q)u_t$. Then
\[
 u_{qc}=\rho_E(q)u_c,\qquad
 u_{cq}=u_c\rho_E(q)
 \qquad(q\in Q,\ c\in C).
\]
Let
\[
 \lambda:\operatorname{End}_{kQ}(E)\longrightarrow
 \operatorname{End}_{kQ}(E)/J(\operatorname{End}_{kQ}(E))\cong k
\]
be the residue map. Since $u_{cd}^{-1}u_cu_d$ is a $kQ$-automorphism of $E$, the scalars
\[
 \alpha(cQ,dQ)=\lambda(u_{cd}^{-1}u_cu_d)^{-1}
\]
form a normalized factor set on $C/Q$. Indeed, the displayed normalization makes this scalar independent of the representatives, and associativity gives the cocycle identity. Write $k_\alpha[C/Q]$ for the twisted group algebra with basis $e_{cQ}$ and multiplication
\[
 e_{cQ}e_{dQ}=\alpha(cQ,dQ)e_{cdQ}.
\]

If $N$ is a $kG$-module with action $\rho_N$, define
\[
 T_c([f])=[\rho_N(c)fu_c^{-1}]
 \qquad(c\in C,\ [f]\in\mathcal M_E(N)).
\]
These operators are well defined, preserve the categorical radical, and satisfy
\[
 T_q=1\quad(q\in Q),\qquad
 T_cT_d=\alpha(cQ,dQ)T_{cd}.
\]
They therefore make $\mathcal M_E(N)$ a $k_\alpha[C/Q]$-module. The factor set is independent of $N$. For a subgroup $D$ of $C$ containing $Q$, the same operators $u_d$ give the restricted factor set $\alpha|_{D/Q}$.

\begin{lemma}\label{lem:multspace}
Let $N$ be an indecomposable $kG$-module with vertex $Q$ and $Q$-source $E$. With the choices above, $\mathcal M_E(N)$ is the Green--Puig multiplicity module of $N$ for $(Q,E)$. It is an indecomposable projective $k_\alpha[C/Q]$-module. Moreover, nonisomorphic indecomposable $kG$-modules with vertex $Q$ and source $E$ have nonisomorphic multiplicity modules over this fixed twisted group algebra.
\end{lemma}

\begin{proof}
Set $A=\operatorname{End}_k(N)$, and write $N\!\downarrow_Q\cong E^m\oplus N_0$ with $E$ absent from $N_0$. The point $\gamma$ of $A^Q=\operatorname{End}_{kQ}(N)$ determined by projection onto a copy of $E$ is local because $E$ has vertex $Q$. Its multiplicity algebra is the corresponding simple quotient of $A^Q$. Hence
\[
 A(Q_\gamma)\cong M_m(k)
 \cong\operatorname{End}_k(\mathcal M_E(N)).
\]
An element of $\N_G(Q)$ stabilizes $\gamma$ exactly when it stabilizes $[E]$, so $\N_G(Q_\gamma)=C$. Under this identification, conjugation on $A(Q_\gamma)$ is implemented by the operators $T_c$ defined above. Thus $\mathcal M_E(N)$ realizes the multiplicity module with the source intertwiners $u_c$ fixing its central extension and factor set; see \cite[2.10 and 6.4]{Puig1988}.

For the reference module $\Ind_Q^C E$, the coset decomposition gives a basis of its multiplicity space on which the operators $T_c$ give the regular $k_\alpha[C/Q]$-module. This is the twisted algebra used in the Green--Puig parametrization. Its indecomposable projective modules parametrize the indecomposable $kG$-modules with vertex $Q$ and source $E$ by \cite[Theorem~1.7 and Remark~1.8]{Harris}; see also \cite[1.3]{Puig}. This proves the remaining assertions.
\end{proof}

The following induction formula is also obtained in \cite[A.6]{Puig}. We give a proof for arbitrary finite groups, since that discussion occurs in an appendix concerning groups of odd order.

\begin{lemma}\label{lem:multind}
Let $H\le G$, let $Q\le H$, and let $E$ be an indecomposable $kQ$-module with vertex $Q$. Set
\[
 C=\Stab_{\N_G(Q)}([E]),\qquad
 D=\Stab_{\N_H(Q)}([E])=C\cap H.
\]
Suppose that $L$ is a simple $kH$-module with vertex $Q$ and source $E$, and that $M=\Ind_H^G L$ is simple with vertex $Q$ and source $E$. Choose the source operators $u_c$ on $C$ as above and use their restrictions on $D$. If $W=\mathcal M_E(L)$ and $U=\mathcal M_E(M)$, then
\[
 U\cong k_\alpha[C/Q]\otimes_{k_{\alpha|_{D/Q}}[D/Q]}W.
\]
\end{lemma}

\begin{proof}
Mackey decomposition gives
\[
 M\!\downarrow_Q\cong
 \bigoplus_{x\in Q\backslash G/H}
 \Ind_{Q\cap{}^xH}^{Q}
 \bigl({}^xL\!\downarrow_{Q\cap{}^xH}\bigr).
\]
If a summand contributes to $\mathcal M_E(M)$, it has a direct summand isomorphic to $E$. Since $E$ has vertex $Q$, relative projectivity forces $Q\le{}^xH$. Also,
\[
 {}^xL\mid\Ind_{{}^xQ}^{{}^xH}({}^xE),
\]
where $\mid$ denotes being isomorphic to a direct summand. Restricting to $Q$ and applying Mackey again shows that $E$ occurs in a term induced from $Q\cap{}^{yx}Q$ for some $y\in{}^xH$. Its vertex forces $Q\cap{}^{yx}Q=Q$. Since these groups have equal order, $Q={}^{yx}Q$, and the term is the indecomposable module ${}^{yx}E$. Consequently,
\[
 E\cong{}^{yx}E.
\]
Writing $y=xhx^{-1}$ with $h\in H$, we obtain $c=yx=xh\in C$ and $cH=xH$.

Conversely, if $c\in C$, then ${}^cL\!\downarrow_Q$ has a direct summand ${}^cE\cong E$. Thus the contributing Mackey double cosets are exactly those represented by elements of $C$. Since $Q\le D$ and $D=C\cap H$, two such representatives give the same double coset precisely when they give the same coset in $C/D$.

In the model $M=kG\otimes_{kH}L$, left multiplication by $C$ permutes the $Q$-orbit summands. Applying the additive functor $\mathcal M_E$ gives
\[
 U=\bigoplus_{cD\in C/D}U_{cD},\qquad U_D=W.
\]
The operators defined above satisfy
\[
 T_c(U_{xD})=U_{cxD}\qquad(c,x\in C),
\]
and, for $d\in D$, the restriction of $T_d$ to $U_D$ is exactly the operator defining the action on $W$, because the same $u_d$ is used. Therefore the map
\[
 k_\alpha[C/Q]\otimes_{k_{\alpha|_{D/Q}}[D/Q]}W
 \longrightarrow U,\qquad
 e_{cQ}\otimes w\longmapsto T_c(w),
\]
is well defined and is an isomorphism on each summand indexed by $C/D$. It is consequently an isomorphism of $k_\alpha[C/Q]$-modules.
\end{proof}

\section{The inducing bound}
For a finite group $G$, a simple $kG$-module $M$ with vertex $Q$, and a subgroup $H$ with $Q\le H\le G$, let
\[
 \mathcal I_M(H\mid Q)=\{[L]:L\text{ is a simple }kH\text{-module with vertex }Q,
 \ \Ind_H^G L\cong M\},
\]
where $[L]$ denotes the isomorphism class of $L$.

\begin{lemma}\label{lem:sourceind}
Let $H\le G$, let $L$ be an indecomposable $kH$-module with vertex $Q$ and $Q$-source $E$, and suppose that $M=\Ind_H^G L$ is indecomposable with vertex $Q$. Then $E$ is, up to $N_G(Q)$-conjugacy, a $Q$-source of $M$.
\end{lemma}

\begin{proof}
Since $L$ is a direct summand of $\Ind_Q^H E$, the module $M$ is a direct summand of $\Ind_Q^G E$. Let $F$ be a $Q$-source of $M$. Then $F$ is a direct summand of
\[
 \Res_Q^G\Ind_Q^G E
 \cong \bigoplus_{x\in Q\backslash G/Q}
 \Ind_{Q\cap{}^xQ}^{Q}({}^xE\!\downarrow_{Q\cap{}^xQ}).
\]
Because $F$ has vertex $Q$, any Mackey term containing $F$ must have $Q\cap{}^xQ=Q$. Hence $x\in N_G(Q)$, and that term is ${}^xE$. Thus $F\cong{}^xE$ for some $x\in N_G(Q)$.
\end{proof}

\begin{theorem}\label{thm:moduleinducing}
Let $G$ be a finite group, let $M$ be a simple $kG$-module with vertex $Q$, and let $Q\le H\le G$. Then
\[
 |\mathcal I_M(H\mid Q)|\le |N_G(Q):N_H(Q)|.
\]
\end{theorem}

\begin{proof}
Set
\[
 A=N_H(Q),\qquad B=N_G(Q).
\]
Choose a $Q$-source $E$ of $M$, and let
\[
 C=N_G(Q,E)=\Stab_B([E]).
\]
For $[L]\in\mathcal I_M(H\mid Q)$, choose a $Q$-source $E_L$ of $L$. By Lemma~\ref{lem:sourceind}, $E_L$ is $B$-conjugate to $E$. Thus the possible $A$-orbits of source classes are indexed by
\[
 C\backslash B/A.
\]

Fix $CbA\in C\backslash B/A$, and let
\[
 F=E^b,\qquad D_b=\Stab_A([F])=A\cap C^b,
\]
and let $\mathcal S_b$ be the set of classes $[L]\in\mathcal I_M(H\mid Q)$ whose source is $A$-conjugate to $F$. If $\mathcal S_b$ is empty, its contribution to the required sum is zero. Suppose that $\mathcal S_b$ is nonempty. After conjugating the source inside $A$, take $F$ as a source for every $L$ in $\mathcal S_b$. Fix one system of source intertwiners for $F$ as in Section~2 and use it for $M$ and every $L$ in $\mathcal S_b$.

Let $U_b$ be the multiplicity module of $M$ relative to $(Q,F)$, and for $[L]\in\mathcal S_b$ let $W_L$ be the multiplicity module of $L$ relative to $(Q,F)$. By Lemma~\ref{lem:multind},
\begin{equation}\label{eq:indmult}
 U_b\cong \Ind_{D_b/Q}^{C^b/Q}W_L.
\end{equation}
The Green--Puig parametrization implies that the $W_L$ are pairwise nonisomorphic indecomposable projective modules as $[L]$ varies in $\mathcal S_b$. Equation \eqref{eq:indmult} gives
\[
 \dim_k W_L=\frac{\dim_k U_b}{[C^b:D_b]}.
\]
The identity double coset in the Mackey decomposition of \eqref{eq:indmult} shows that each $W_L$ is a direct summand of $\Res_{D_b/Q}^{C^b/Q}U_b$. By the Krull--Schmidt theorem, a decomposition of this restriction into indecomposable modules contains at least one representative of every isomorphism class $[W_L]$. Consequently,
\[
 |\mathcal S_b|\frac{\dim_k U_b}{[C^b:D_b]}
 =\sum_{[L]\in\mathcal S_b}\dim_k W_L
 \le \dim_k U_b,
\]
and hence
\begin{equation}\label{eq:onecoset}
 |\mathcal S_b|\le [C^b:A\cap C^b].
\end{equation}

Summing \eqref{eq:onecoset} over $C\backslash B/A$ gives
\[
 |\mathcal I_M(H\mid Q)|
 \le \sum_{CbA\in C\backslash B/A}[C^b:A\cap C^b]
 =[B:A]
 =|N_G(Q):N_H(Q)|,
\]
because $CbA$ contains exactly $[C^b:A\cap C^b]$ left $A$-cosets.
\end{proof}

\begin{corollary}\label{thm:inducing}
Let $G$ be a finite $p$-solvable group, let $\varphi\in\IBr_p(G)$ have vertex $Q$, and let $Q\le H\le G$. Then
\[
 |I_\varphi(H\mid Q)|\le |N_G(Q):N_H(Q)|.
\]
\end{corollary}

\begin{proof}
Let $M$ be the simple $kG$-module affording $\varphi$. If $\eta\in I_\varphi(H\mid Q)$ and $L_\eta$ is the corresponding simple $kH$-module, then $\eta^G=\varphi$. The Brauer character of $\Ind_H^G L_\eta$ is therefore irreducible. Since the Brauer character of a module is the sum of the characters of its composition factors, this induced module is simple. Hence
\[
 \Ind_H^G L_\eta\cong M.
\]
Thus $\eta\mapsto[L_\eta]$ embeds $I_\varphi(H\mid Q)$ in $\mathcal I_M(H\mid Q)$, and Theorem~\ref{thm:moduleinducing} gives the result.
\end{proof}

\section{Lifts for odd primes}
Let $G$ be a finite $p$-solvable group, and set $\pi=p'$. Then $G$ is $\pi$-separable and $I_\pi(G)=\IBr_p(G)$. We use Navarro vertices as defined in \cite[Section~2]{ZCJW}. For $\varphi\in\IBr_p(G)$, a $p$-subgroup $Q$, and $\delta\in\Irr(Q)$, let
\[
 L_\varphi(Q,\delta)=\{\chi\in L_\varphi:(Q,\delta)\text{ is a Navarro vertex of }\chi\},
\]
and let $\N_G(Q,\delta)$ be the stabilizer of $\delta$ in $\N_G(Q)$.

Theorem~B of \cite{ZCJW} reduces the bound for a fixed Navarro vertex to inducing inequalities in subgroups. It applies to an arbitrary set of primes $\pi$.

\begin{theorem}\label{thm:fixedpair}
Let $p$ be a prime, let $G$ be a finite $p$-solvable group, and let $\varphi\in\IBr_p(G)$ have vertex $Q$. If $\delta\in\Irr(Q)$, then
\[
 |L_\varphi(Q,\delta)|\le |\N_G(Q):\N_G(Q,\delta)|.
\]
\end{theorem}

\begin{proof}
Let $Q\le R\le S\le G$ and let $\sigma\in I_\varphi(S\mid Q)$. Then $\sigma$ has vertex $Q$. Since subgroups of a $p$-solvable group are $p$-solvable, Corollary~\ref{thm:inducing} applied to $R\le S$ gives
\[
 |I_\sigma(R\mid Q)|\le |\N_S(Q):\N_R(Q)|.
\]
Thus the hypothesis of \cite[Theorem~B]{ZCJW} holds for every such pair $R\le S$, and that theorem gives the required inequality.
\end{proof}

\begin{theorem}\label{thm:cossey}
Let $p>2$, let $G$ be a finite $p$-solvable group, and let $\varphi\in\IBr_p(G)$ have vertex $Q$. Then
\[
 |L_\varphi|\le |Q:Q'|.
\]
\end{theorem}

\begin{proof}
For $p>2$, the generalized vertex pairs of every lift of $\varphi$ are conjugate by \cite[Theorem~1.2]{CLvertex}, and the corresponding vertex character is linear \cite[Lemma~4.3]{CLvertex}. In the Brauer-character setting $\pi=p'$, so $2\in\pi$; in this case the Navarro-vertex notation $L_\varphi(Q,\delta)$ is compatible with the Cossey--Lewis vertex-pair notation, as observed in \cite[Introduction]{ZCJW}. Let $\delta_1,\ldots,\delta_t$ represent the $\N_G(Q)$-orbits on $\Lin(Q)$. Then
\[
 L_\varphi=\coprod_{i=1}^t L_\varphi(Q,\delta_i).
\]
By Theorem~\ref{thm:fixedpair},
\[
 |L_\varphi|
 \le \sum_{i=1}^t |\N_G(Q):\N_G(Q,\delta_i)|
 =|\Lin(Q)|
 =|Q:Q'|.\qedhere
\]
\end{proof}

\section*{Acknowledgements}
Yang was supported by Simons Foundation Award \#918096.

The author used ChatGPT to assist with developing and checking the proofs and revising the exposition. The author is responsible for the mathematical content.

\section*{Disclosure Statement}
The author reports no conflict of interest.

\section*{Data Availability}
No data were used in this work.

\end{document}